\documentclass[11pt]{article}
\usepackage{amsmath,amsfonts,amsthm,epsf,stackrel}
\usepackage{graphicx,subfigure}
\usepackage{natbib}

\RequirePackage{natbib}

\usepackage[margin=2cm]{geometry}

\newcommand{\ignore}[1]{}
\newcommand{\tendp}{\stackrel{p}{\longrightarrow}}

\newcommand{\ve}{\varepsilon}

\newtheorem{theorem}{Theorem}[section]
\newtheorem{corollary}{Corollary}[section]

\newtheorem{lemma}{Lemma}[section]

\begin{document}

\title{The conditionality principle in high-dimensional regression}
\author{
	David Azriel
}
\vspace{0.5in}
\date{\today}
\maketitle

\begin{abstract}	
Consider a high-dimensional linear regression problem, where the number of covariates is larger than the number of observations and the interest is in estimating the conditional variance of the response variable given the covariates. A conditional and unconditioned framework are considered, where conditioning is with respect to the covariates, which are ancillary to the parameter of interest. In recent papers, a consistent estimator was developed in the unconditional framework when the marginal distribution of the covariates is normal with known mean and variance. In the present work, a certain Bayesian hypothesis test is formulated under the conditional framework, and it is shown that the Bayes risk is a constant. This implies that no consistent estimator exists in the conditional framework. However, when the marginal distribution of the covariates is normal, the conditional error of the above consistent estimator converges to zero, with probability converging to one. 
It follows that even in the conditional setting, information about the marginal distribution of an ancillary statistic may have a significant impact on statistical inference. The practical implication in the context of high-dimensional regression models is that additional observations, where only the covariates are given, are potentially very useful and should not be ignored. This finding is most relevant to semi-supervised learning problems where covariate information is easy to obtain.
\end{abstract}

\section{Introduction}

An ancillary statistic is one whose distribution does not depend on the parameters of the model. 
In \citet{Cox_1974} notation and words, if $C$ is an ancillary statistic then ``the conditionality principle is that the conclusion about the parameter of interest is to be drawn as if $C$ were fixed at its observed value $c$'' (p. 38). This principle has two implications:
\begin{enumerate}
	\item {Conditional inference:} Statistical inference should be conditioned on an ancillary statistic. 
	\item {Ignorability of the marginal distribution:} The true marginal distribution of an ancillary statistic should be ignored in any estimation procedure.
\end{enumerate}
Focusing on the conditional inference implication, \citet{Brown} presents an ancillarity paradox, where in a certain regression problem, the standard estimator is admissible in the conditional setting for every value of the ancillary statistic but it is inadmissible in the unconditional setting. 
Brown argues that the common practice to consider only conditional inference is sometimes misleading and the unconditional risk function should also be accounted for. 
%See \citet{Ghosh2010} for further discussions.
%He concludes that ``the admissibility results of this paper contradict the widely held notion that statistical inference in the presence of ancillary statistic should be independent of those ancillary statistic''. 
The present work goes one step further with respect to the above ignorability implication. It is shown, under the framework of the conditional setting, that an estimator using the true marginal distribution of an ancillary statistic has a vanishing conditional error with high probability. 
%Notice that in the conditional setting the Bayes risk is random. 
That is, even if one carries out conditional inference, one can still benefit from the marginal distribution of an ancillary statistic, rejecting the ignorability implication. Moreover, while in Brown's setting the advantage of the new estimator becomes negligible as the sample size grows; in the high-dimensional asymptotic regime considered here, the conditional error goes to zero with high probability. This is because in high-dimensional regression the sampling distribution itself changes with sample size.

The regression linear model is
%\begin{equation}\label{eq:lin_model}
${ Y}={ X} { \beta} + { \ve}$,
%\end{equation}
where ${ Y}=(y_1,\ldots,y_n)^T$ is a column vector of the response variables, ${ X}$ is an $n \times p$ matrix of covariates, ${ \beta} = (\beta_1,\ldots,\beta_p)^T$ is a vector of unknown parameters and ${ \ve}=(\ve_1,\ldots,\ve_n)^T$ is a vector of the residuals. It is assumed that $(y_1,{ x}_1), \ldots, (y_n,{ x}_n)$ are independently and identically distributed, where ${ x}_i$ is the $i$-th row of the matrix ${ X}$; $({ x},y)$ denotes a generic observation that should be read as
$({ x}_i,y_i)$ for some $i$. Here the focus is on estimating $\sigma^2={\rm var}(y|{ x})$.
When the linear model is true, i.e., ${\rm E}(y|{ x})={ x}^T{ \beta}$, then ${ \beta}$ is the coefficient of the conditional expectation; when the conditional expectation is not linear in ${ x}$, then ${ \beta}$ is still meaningful in the sense it can be defined as the coefficient of the best linear predictor; see \citet{Buja} for exact definitions.

Consider an additional sample of size $m$ where only observations from the marginal distribution of ${ x}$ are given. In the machine learning literature, this situation is called semi-supervised learning; see \citet{Zhou2014} and references therein. Such a situation arises when the $y$ observations are costly and the ${ x}$'s are easy to obtain. A typical example is web document classification, where the classification is done by a human agent while there are many more unlabeled online documents. In such situations, $m$ is much larger than $n$, and hence the marginal distribution of ${ x}$ can be assumed known.

While there is a large body of literature on this problem in the machine learning world, there are very few statistical papers concerning this topic. The standard approach in the statistical literature may be best summarized by the following quote from \citet{Little1992}:
\begin{quote}
	The related problem of missing values in the outcome $Y$ was prominent in the early history of missing-data methods, but is less interesting in the following sense: If the $X$'s are complete and the missing values of $Y$ are missing at random, then the incomplete cases contribute no information to the regression of $Y$ on $X_1, \ldots, X_p$.
\end{quote}
This approach is justified by the conditionality principle since ${ X}$ is ancillary to the parameters of interest. More generally, in the context of causal inference, 
\citet{Janzing} argue that knowing the marginal distribution of ${ X}$ is useless for the study of the conditional distribution of ${ Y}$ given ${ X}$ when ${ X}$ has a causal effect on ${ Y}$. 
On the other hand, \citet{Buja} show that ${ X}$ is ancillary with respect to ${ \beta}$ if and only if the linear model is actually true. Indeed, for a low-dimensional regression model, \citet{Chakrabortty} and \citet{Azriel1} construct a semi-supervised estimator that asymptotically dominates the least squares estimator as the latter is based only on the labeled data set. The asymptotic variance of the new estimators are smaller than that of the least squares estimator when a certain non-linearity condition holds; otherwise, the new estimates are equivalent to the least squares estimator. In other words, improvement can be made only in the non-linear case where ${ X}$ is no longer ancillary. 

In high-dimensional regression, $p$ is larger than $n$, i.e., there are more parameters than observations. In this setting, the lasso estimator suggested by \citet{Tibshirani1996} has gained much popularity and many extensions were suggested. That line of research is related to sparsity assumptions where most of the parameters are assumed to be zero or close to zero. When those assumptions do not hold, then estimation of the entire vector of ${ \beta}$ is not feasible. However, it is possible to estimate the signal-to-noise-ratio and $\sigma^2$. \citet{Dicker}, \citet{Dicker1} and \citet{Janson} suggest estimators when assuming that the marginal distribution of ${ x}$ is standard normal with independent entries. In the context of semi-supervised learning and when the marginal distribution of ${ x}$ is $N({ \mu},{ \Sigma})$ say, then by a linear transformation, which does not change $\sigma^2$ or the signal-to-noise ratio, ${ x}$ could become standard normal with independent entries, justifying their assumptions. 

The estimator of \citet{Dicker} is based on moments, as described below in Section \ref{sec:prel}; while \citet{Dicker1} suggest a maximum likelihood estimator. \citet{Janson} follow a different approach and study an optimization problem for which the resulting estimator is unbiased and has small variance.   

The interesting point is that in these three works the marginal distribution of ${ x}$ plays an important role in the estimation procedure even though it is ancillary to the parameters of interest. This point is the motivation for the present paper. The focus here is on the estimator of \citet{Dicker} since it is easier to analyze, although it might be inferior to the estimators of \citet{Dicker1} and \citet{Janson}. It is shown that under an asymptotic regime where $p/n \to c \in [1,\infty)$, the conditional error of Dicker's estimate is arbitrarily close to zero with probability converging to one. On the other hand, in the conditional setting when $x$ is assumed constant, no consistent estimator exists. 
This result  implies that information on the marginal distribution of an ancillary statistic can be helpful even in the conditional setting contradicting the ignorability implication of the conditionality principle that was mentioned above.

\section{Preliminaries}\label{sec:prel}

%Let $(y_1,x_1),...,(y_n,x_n)$ be i.i.d where ${ x}$ is a (row) vector of covariates in ${\mathbb R}^p$ and $y\in \mathbb{R}$ is the response; $({ x},y)$ denotes a ``generic'' observation that should be read as
%$({ x}_i,y_i)$ for some $i$. 

Assume that $(y_1,x_1),...,(y_n,x_n)$ are sampled from the distribution
\begin{equation}\label{eq:model}
y|{ x}={ x}_0 \sim N({ \beta}^T { x}_0,\sigma^2),~{ x} \sim N({ 0}, { I}),
\end{equation}
for some unknown vector ${ \beta} \in \mathbb{R}^p$ and parameter $\sigma^2$, and where ${ I}$ is the $p\times p$ identity matrix.
The purpose is to estimate $\sigma^2$.  
The present work studies model \eqref{eq:model} under an asymptotic regime in which both the dimension $p=p_n$ and the number of observations $n$ converge to infinity with $p/n \to c$ as $n \to \infty$, where $c \in [1,\infty)$. It is also assumed, as in \citet{Dicker}, that $||{ \beta}||^2$ and  $\sigma^2$ are of order of a constant as $n \to \infty$. Since this constant is unknown and no sparsity assumptions are made,  the parameter space for $\beta$ is $\mathbb{R}^p$ and for $\sigma^2$ is $\mathbb{R}^+$.
% The latter assumption ensures that the signal-to-noise ratio is non-trivial. 

Estimation of $\sigma^2$ is considered in two settings:
\begin{itemize}
	\item Conditional setting: the distribution is conditioned upon the observations $x_1,\ldots,x_n$, or equivalently upon the matrix $X$. That is, $y_i \sim N(\mu_i,\sigma^2)$, with $\mu_i = {\rm E}(y_i|x_i)$, $(i=1,\ldots,n)$. Since $p\ge n$, there are no constraints on the $\mu$'s, i.e., the parameter space for the $\mu$'s is $\mathbb{R}^n$, as ${ X}$ is of full rank with probability 1. To clarify the notation we use ${\rm pr}_X$, ${\rm E}_{X}$ and ${\rm var}_{X}$ to denote conditional probability, expectation and variance; that is, for a given matrix $X$, these three quantities are numbers and not random variables.
	\item Unconditional setting as above. Here one can use the known marginal distribution of the ${ x}$'s. 
\end{itemize}

The parameter space for ${\beta}$ is $\mathbb{R}^p$ and therefore, in the conditional setting, since $p \ge n$ the $\mu$'s are unrestricted and ${ X}$ plays no role. This is different from low-dimensional regression, where conditioning on ${ X}$ restricts the $\mu$'s to lie in a $p$-dimensional sub-space of $\mathbb{R}^n$. This is the reason why there is such a difference in conditioning in low- and high-dimensional regression. Unlike the unconditional setting, in the conditional one, it is shown below that there exists no consistent estimator for $\sigma^2$. This implies that in a high-dimensional regression model, when estimation of $\sigma^2$ is of interest, the unconditional setting should be preferred. It is important to notice that no sparsity assumptions on ${ \beta}$ are imposed. If such assumptions are made, then the $\mu$'s belong to a restricted set of $\mathbb{R}^p$ and consistent estimation is possible also in the conditional setting; see, e.g., \citet{Sun}. 

The parameter $\sigma^2$ is identifiable in both the conditional and unconditional settings, since different values of $\sigma^2$ correspond to different probability measures in both. On the other hand, the vector ${ \beta}$ is identifiable only in the unconditional setting. Therefore, a better comparison can be made with respect to $\sigma^2$ and this is the focus here.

\citet{Dicker} considered the unconditional setting and noticed that under model \eqref{eq:model},
\[
{\rm E}\left(\frac{1}{n}||{ Y}||^2\right)=||{ \beta}||^2+\sigma^2~~~,~~~{\rm E}\left(\frac{1}{n^2}||{ X}^T { Y}||^2\right)=\frac{p+n+1}{n}||{ \beta}||^2 + \frac{p}{n} \sigma^2,
\]
where these expectations are unconditional.
This leads to the estimator 
\begin{equation}\label{eq:Dicker}
\hat{\sigma}^2_{Dicker}=\frac{p+n+1}{n(n+1)}||{ Y}||^2-\frac{1}{n(n+1)}||{ X}^T{ Y}||^2.
\end{equation}
%Under the unconditional setting a consistent estimator was developed by \citet{Dicker} (see also \citet{Janson}).
%For the conditional setting,  a Bayesian hypothesis testing and is considered and it is shown that the corresponding Bayes risk is constant. However, a rule based on Dicker's estimate attains a smaller Bayes risk. This is surprising since ${ x}$ is ancillary in the above model and so the information about the marginal distribution of ${ x}$ should not help us in estimating $\sigma^2$.
%Dicker's estimate is
%\begin{equation}\label{eq:Dicker}
%\hat{\sigma}^2_{Dicker}=\frac{p+n+1}{n(n+1)}||{ Y}||^2-\frac{1}{n(n+1)}||{ X}^T{ Y}||^2,
%\end{equation}
%Here the focus is on $\hat{\sigma}^2_{Dicker}$.
It follows that ${\rm E}(\hat{\sigma}^2_{Dicker})=\sigma^2$; the variance is 
\begin{equation}\label{eq:V_dicker}
{\rm var}(\hat{\sigma}^2_{Dicker})=\frac{2}{n}\left\{\frac{p}{n}(\sigma^2+||{ \beta}||^2)^2+\sigma^4+||{ \beta}||^4\right\}\{1+O(1/n)\}.
\end{equation}
Conditional analysis of Dicker's estimate is carried out below. 

\section{The Bayes risk}

Consider now the conditional setting. Suppose we want to determine whether $\sigma^2=\sigma_0^2$ or $\sigma^2=\sigma_1^2$. Let $\delta: { Y} \mapsto \{0,1\}$ be a rule and denote the true value by $J \in \{0,1\}$. The risk is $R_\delta={\rm pr}_X(\delta \ne J)$, and it is a function of the unknown parameters $\mu_1,\ldots,\mu_n,\sigma^2$. 
We show below that for any $\delta$,
\begin{equation}\label{eq:int}
\int R_\delta d\pi(\mu_1,\ldots,\mu_n,\sigma^2) \ge 1/2,
\end{equation}
where $\pi$ is a certain probability measure on the parameter space $\mathbb{R}^n \times \{\sigma_0^2,\sigma_1^2\}$ to be defined below.
The probability measure $\pi$ in \eqref{eq:int} can be thought of as a prior density for the parameters. For any rule $\delta$, the integral in \eqref{eq:int} is bounded below by the Bayes risk. However, in the next section we show that under a rule based on Dicker's estimate, the integral in \eqref{eq:int} can be arbitrarily small. This rule uses the knowledge of the marginal distribution of ${ X}$. This demonstrates the usefulnesses of the information on the distribution of an ancillary statistic in this context.

In order to show \eqref{eq:int}, it is enough to compute the Bayes risk. The idea is to define $\pi$ as a mixture of two priors 
$\mu_1,\ldots,\mu_n \sim N(0,\eta_0^2), \sigma^2=\sigma_0^2$ and $\mu_1,\ldots,\mu_n \sim N(0,\eta_1^2), \sigma^2=\sigma_1^2$,    
where the notation $\sim$ means that the random variables are independent and identically distributed.
When $\eta_0^2+\sigma_0^2=\eta_1^2+\sigma_1^2$, the $y$'s have the same marginal distributions under the two priors and therefore no Bayesian procedure can distinguish between them, leading to \eqref{eq:int}.  

Specifically, consider a Bayesian formulation where $J \sim {\rm Bernoulli}(1/2)$. Given $J=0$,  
\[
\mu_1,\ldots,\mu_n \sim N(0,\eta_0^2), ~\sigma^2=\sigma_0^2;
\]
and similarly for $J=1$ with $\eta_1^2,\sigma_1^2$ replacing $\eta_0^2,\sigma_0^2$. Let $\pi$ be the induced probability measure on $\mu_1,\ldots,\mu_n,\sigma^2$. Lemma \ref{lem} computes the Bayes rule under $\pi$. Specifically, it is shown that if $\eta_0^2+\sigma_0^2=\eta_1^2+\sigma_1^2$ the posterior probability of $J$ is the same as the prior, implying \eqref{eq:int}. 
\begin{lemma}\label{lem}
	If $\eta_0^2+\sigma_0^2=\eta_1^2+\sigma_1^2$ then
	\[
	{\rm pr}_{X,\pi}(J=0 \mid y_1,\ldots,y_n)={\rm pr}_{X,\pi}(J=1 \mid y_1,\ldots,y_n)=1/2,	
	\]
	where ${\rm pr}_{X,\pi}$ denotes the conditional probability under the prior $\pi$.
	%Otherwise, $P_{x,\pi}(J=0|y_1,\ldots,y_n)$ is higher if
	%\[
	%\frac{1}{n}\sum_{i=1}^n y_i^2 > \frac{\log\left(\frac{\eta_0^2+\sigma_0^2}{\eta_1^2+\sigma_1^2}\right)}{\frac{1}{\eta_1^2+\sigma_1^2}-\frac{1}{\eta_0^2+\sigma_0^2}}.
	%\]
\end{lemma}

Since \eqref{eq:int} holds for any rule $\delta$, it follows that there exists no consistent estimator for $\sigma^2$ in the conditional setting. This result is summarized in the next corollary.
\begin{corollary}\label{cor}
	Under model \eqref{eq:model} in the conditional setting there exists no consistent estimator for $\sigma^2$.
\end{corollary}

\section{Dicker's estimate in the conditional setting}

Recall Dicker's estimate, which is defined in \eqref{eq:Dicker} and consider the problem of the previous section to determine whether $\sigma^2=\sigma_0^2$ or $\sigma^2=\sigma_1^2$. Assume without loss of generality that $\sigma_1^2 > \sigma_0^2$. Define a rule based on Dicker's estimate
\[
\delta_{Dicker} = \left\{ \begin{array}{cc} 1 & \hat{\sigma}^2_{Dicker} > \frac{\sigma_0^2 +\sigma_1^2}{2} \\ 0 & \hat{\sigma}^2_{Dicker} \le  \frac{\sigma_0^2 +\sigma_1^2}{2} \end{array} \right. .
\] 
Define the conditional error $R_{\delta_{Dicker}} = {\rm pr}_X( \delta_{Dicker} \ne J)$. We show below that with high probability (over ${ X}$), $R_{\delta_{Dicker}}$ is small. 
%{\bf DD DELETE} ``This implies that the Bayes risk of $\delta_{Dicker}$ is smaller than that of the Bayes rule with high probability.'' 

\begin{theorem}\label{thm}
	Consider model \eqref{eq:model}, and assume  that $p/n \to c$ as $n \to \infty$ for $c \in [1,\infty)$, and that $\sigma^2,||{ \beta}||^2$ are bounded as $n \to \infty$. Then, there exist a sequence of sets $A_n$ and a constant $C$, such that ${\rm pr}({ X} \in A_n) \to 1$ as $n \to \infty$ and for ${ X} \in A_n$ and any $\xi > 0$,
	\begin{equation}\label{eq:thm}
	{\rm pr}_X\left(|\hat{\sigma}^2_{Dicker}-\sigma^2| \ge \xi\right) \le  C \frac{g(c,n,{ \mu},\sigma^2)}{\xi^2 \sqrt{n}},
	\end{equation} 
	where $g(c,n,{ \mu},\sigma^2)=1+ 2(c+1)\left\{\left(\frac{1}{n}||{ \mu}||^2\right)^2+\sigma^4\right\}+\frac{4 \sigma^2}{n}||{ \mu}||^2 + 2\sigma^4$.
\end{theorem}	

Consider the conditional probability of the event that $\delta_{Dicker}=1$ but $\sigma^2_0$ is actually true, i,e., $\sigma^2=\sigma_0^2$; then,
\[
{\rm pr}_X\{\hat{\sigma}^2_{Dicker} \ge  (\sigma^2_1+\sigma^2)/2 \} \le {\rm pr}_X\left(|\hat{\sigma}^2_{Dicker}-\sigma^2| \ge \xi \right),  
\]
where $\xi=(\sigma^2_1+\sigma^2)/2$. 
%{\bf DD Delete} ``Theorem \ref{thm} implies that for ${ X}\in A_n$,
%\[
%\int {\rm pr}_X\{\hat{\sigma}^2_{Dicker}  \ge  (\sigma^2_1+\sigma^2)/2 \} d\pi(\mu_1,\ldots,\mu_n,\sigma^2) \le C\int g(c,n,{ \mu},\sigma^2)d\pi(\mu_1,\ldots,\mu_n,\sigma^2)/(\xi^2\sqrt{n}).
%\] 
%Since $g(c,n,{ \mu},\sigma^2)$ has finite expectation with respect to $\pi$, the conditional probability converges to zero.'' 
It follows that for ${ X}\in A_n$, $R_{\delta_{Dicker}}$ converges to zero at rate $\sqrt{n}$ provided that $||\mu||^2$ is bounded. However, Theorem \ref{thm} cannot be used to calculate the Bayes risk of $\delta_{Dicker}$ since under the prior $\pi$, $||\beta||^2$ is unbounded with high probability and the conditions of the theorem are not satisfied.   

\section{Discussion}

%To sum up, a Bayesian hypothesis testing is considered in the conditional setting, where 
%it is shown that the corresponding Bayes risk is constant, while a rule based on Dicker's estimate attains a smaller Bayes risk with high probability. More generally, by Corollary \ref{cor}, there exists no consistent estimator for $\sigma^2$ in the conditional setting, whereas Dicker's estimate is consistent in the unconditional setting. 
%Since ${ x}$ is ancillary in the above model, this finding demonstrates that the ignorability implication of the conditionality principle that was mentioned in the Introduction is inadequate under model \eqref{eq:model}, at least in the context of the estimation of $\sigma^2$. 
The argument in the previous section is valid in an asymptotic regime where $n \to \infty$. In Section \ref{sec:sim} below, $R_{\delta_{Dicker}}$ is evaluated for certain values of the parameters using simulations. It is demonstrated that even for a relatively small sample size, such as $n=100$,  the probability of $R_{\delta_{Dicker}}$ being smaller than 0.5 is practically one. 

Unlike in \citet{Brown}, here it is shown that an estimator based on information about the marginal distribution of an ancillary statistic is helpful also in the conditional setting. On the other hand, Brown considers the conditional setting for every possible value of ${ X}$, while here the claims are given only with high probability with respect to the distribution of ${ X}$.  

The results here are close in spirit to the work of \citet{Robins1997}. They consider a certain semi-parametric problem and show that when the marginal distribution of an ancillary statistic is known, a consistent estimator exists, but otherwise consistent estimation is impossible. 
The current work shows a similar phenomenon in the context of the widely-used and much simpler high-dimensional regression model. 
%{\bf DD DELETE} ``the current work goes one step further by presenting a paradoxical situation: one can construct a procedure that is based on the marginal distribution of an ancillary statistic and has a smaller Bayes risk than the Bayes rule''. 
As pointed out by
\citet{Robins2000}, in reference to the work of \citet{Robins1997}, these results call for rethinking fundamental principles in the presence of modern data sets where high-dimensional or infinite-dimensional models are natural.  

%The present work is another example where low-dimensional models are different than high-dimensional models with respect to the conditionality principle.

The purpose of the present work is not only to make general comments on the conditionality principle but also to emphasize the importance of the semi-supervised framework in the context of high-dimensional regression. In many data sets, one can easily obtain observations where only the covariates are given but labeled observations are costly. For example, tracking  Parkinson's disease progression requires time-consuming physical 
examinations by trained medical staff; however, self-administered speech tests are easy to run and can be used to predict the symptom progression \citep{Tsanas}. In a different context, counting the number of people who are homeless in census tracts is costly but one can use publicly available data to infer this information in tracts where it is unknown \citep{Kriegler}. Other examples for semi-supervised classification can be found in \citet{ALMA}. In such situations, where there are many more unlabeled than labeled observations, one can use the unlabeled data to estimate the marginal distribution of the covariates and to improve inference. The argument presented here emphasizes the importance of the marginal distribution of the ${x}$'s for the estimation of $\sigma^2$ in high-dimensional regression. As mentioned in the introduction, semi-supervised regression problems have received little attention in the statistical literature. The current work aims at changing this situation.

\section{Proofs}

\subsection*{Proof of Lemma 1}
All the probabilities and expectations below are with respect to $\pi$ and the subscript $\pi$ is suppressed in the notation. The difference of the log posterior probabilities is 
\begin{multline*}
\log({\rm pr}_X(J=0\mid y_1,\ldots,y_n))-\log({\rm pr}_X(J=1\mid y_1,\ldots,y_n))\\
=\log(f(y_1,\ldots,y_n\mid J=0))-\log(f(y_1,\ldots,y_n\mid J=1)),
\end{multline*}
where $f$ is the conditional density; notice that the prior does not play a role since ${\rm pr}_X(J=0)={\rm pr}_X(J=1)$. Furthermore,
\[
f(y_1,\ldots,y_n\mid J=0)={\rm E}\{ f(y_1,\ldots,y_n\mid { \mu},\sigma^2) \mid  J=0\},
\]
where the expectation is over ${ \mu}$ (the $y$'s are constants and $\sigma^2$ is also constant given $J=0$). Now,
\[
{\rm E}_X\{ f(y_1,\ldots,y_n\mid { \mu},\sigma^2) \mid  J=0\}={\rm E}\left[ \left({2 \pi \sigma_0^2}\right)^{-n/2} \exp\left\{-\sum_{i=1}^n (y_i - \mu_i)^2/2\sigma_0^2\right\} \right],
\]
where the latter expectation is over ${ \mu} \sim N(0,\eta_0^2 { I})$, i.e., the $y$'s are constants.
To compute the latter expectation, consider the integral 
\[
\int_{-\infty}^{\infty}\left({2 \pi \sigma_0^2}\right)^{-1/2}\exp\left\{-\frac{(y-\mu)^2}{2 \sigma_0^2}\right\} \left(2 \pi \eta_0^2\right)^{-1/2}\exp\left(-\frac{\mu^2}{2 \eta_0^2}\right) d\mu.
\]
By trivial algebra we have
\[
\frac{(y-\mu)^2}{\sigma_0^2}+\frac{\mu^2}{\eta_0^2}=\left(\frac{1}{\sigma_0^2}+\frac{1}{\eta_0^2}\right)\left(\mu - \frac{y}{1+\sigma_0^2/\eta_0^2}\right)^2 + \frac{y^2}{\eta_0^2 + \sigma_0^2}.
\]
Therefore,
\[
\int_{-\infty}^{\infty}\left(2 \pi \sigma_0^2\right)^{-1/2}\exp\left\{-\frac{(y-\mu)^2}{2 \sigma_0^2}\right\} \left(2 \pi \eta_0^2\right)^{-1/2}\exp\left(-\frac{\mu^2}{2 \eta_0^2}\right) d\mu
=\frac{1}{\sqrt{2\pi(\sigma_0^2+\eta_0^2)}} \exp\left(-\frac{y^2}{\eta_0^2+\sigma_0^2}\right).
\]
Hence, 
\[
\log\{f(y_1,\ldots,y_n\mid J=0)\}=-\frac{n}{2}\log(2\pi) -\frac{n}{2} \log(\eta_0^2+\sigma_0^2)  - \frac{\sum_{i=1}^n y_i^2}{2(\eta_0^2+\sigma_0^2)},
\]
and
\begin{multline*}
\log\{f(y_1,\ldots,y_n\mid J=0)\}-\log\{f(y_1,\ldots,y_n\mid J=1)\}\\
=\frac{n}{2} \log\left(\frac{\eta_1^2+\sigma_1^2}{\eta_0^2+\sigma_0^2}\right)+\frac{1}{2} \sum_{i=1}^n y_i^2\left(\frac{1}{\eta_1^2+\sigma_1^2}-\frac{1}{\eta_0^2+\sigma_0^2}\right).
\end{multline*}
Therefore, when $\sigma_1^2 + \eta_1^2=\eta_0^2 + \sigma_0^2$, the posterior probability that $J=0$ and $J=1$ are equal no matter what the values of the $y$'s are. Otherwise,
the posterior probability that $J=0$ is higher when
\[
\frac{1}{n}\sum_{i=1}^n y_i^2 \ge \frac{\log\left(\frac{\eta_0^2+\sigma_0^2}{\eta_1^2+\sigma_1^2}\right)}{\frac{1}{\eta_1^2+\sigma_1^2}-\frac{1}{\eta_0^2+\sigma_0^2}}.
\]
%\qed

\subsection*{Proof of Corollary 1}
Assume that a consistent estimate exists and denote it by $\hat{\sigma}_*^2$. 
Without loss of generality, assume that $\sigma_1^2 > \sigma_0^2$; now, define the corresponding rule 
\[
\delta_{*} = \left\{ \begin{array}{cc} 1 & \hat{\sigma}^2_{*} > \frac{\sigma_0^2 +\sigma_1^2}{2} \\ 0 & \hat{\sigma}^2_{*} \le \frac{\sigma_0^2 +\sigma_1^2}{2} \end{array} \right. .
\] 
By consistency, $R_{\delta_{*}}$ converges to zero for any sequence $\{\mu_i\}_{i=1}^\infty$ and for $\sigma^2 \in \{\sigma_0^2,\sigma_1^2\}$. However, since $R_{\delta_{*}}$ is bounded, by Lebesgue's dominated convergence theorem, 
\[
\lim_{n \to \infty} \int R_{\delta_{*}} d\pi(\mu_1,\ldots,\mu_n,\sigma^2) =0,
\]
contradicting (4).%\eqref{eq:int}. %\qed 

\subsection*{Proof of Theorem 1}
The proof is based on the decomposition
\begin{equation}\label{eq:E_Var}
{\rm var}(\hat{\sigma}^2_{Dicker})= E \{{\rm var}_X(\hat{\sigma}^2_{Dicker})\}+ {\rm var} \{{\rm E}_X (\hat{\sigma}^2_{Dicker})\},
\end{equation}
which connects the unconditional variance of $\hat{\sigma}^2_{Dicker}$ and the conditional expectation and variance. Also, by (3), there exists a constant $C_1$ such that for every $n$,
\begin{equation}\label{eq:V_bound}
{\rm var}(\hat{\sigma}^2_{Dicker}) \le C_1 V/n,
\end{equation} 
where $V=2\left\{c(\sigma^2+||{ \beta}||^2)^2+\sigma^4+||{ \beta}||^4\right\}$; notice that under our asymptotic regime, $V$ is bounded in $n$. 

%We have that
%\begin{equation}\label{eq:E_var}
%{\rm var}(\hat{\sigma}^2_{Dicker})= E \{var (\hat{\sigma}^2_{Dicker}|{ X})\}+ var \{E (\hat{\sigma}^2_{Dicker}|{ X})\}.
%\end{equation}
%Since both summands are positive, then when ${\rm var}(\hat{\sigma}^2_{Dicker})$ converges to zero, then also these summands converge. Now, since 
%\begin{equation*}
%{\rm E}\{E (\hat{\sigma}^2_{Dicker}|{ X})\}=E (\hat{\sigma}^2_{Dicker})=\sigma^2 \text{ and }
%{\rm var}\{E (\hat{\sigma}^2_{Dicker}|{ X})\} \to 0,
%\end{equation*}
%then $E (\hat{\sigma}^2_{Dicker}|{ X}) \tendp \sigma^2$. Also, by Markov's inequality,
%\[
%P\{var (\hat{\sigma}^2_{Dicker}|{ X}) \ge  \ve \} \le \frac{{\rm E}\{var (\hat{\sigma}^2_{Dicker}|{ X})\} }{\ve} \to 0,
%\]
%and hence $var (\hat{\sigma}^2_{Dicker}|{ X}) \tendp 0$. The conditional Markov's inequality implies for $\delta > 0$,
%\[
%P\left(|\hat{\sigma}^2_{Dicker}-\sigma^2| \ge \delta |{ X}\right) \le \frac{{\rm E}\left\{\left(\hat{\sigma}^2_{Dicker}-\sigma^2\right)^2 |{ X}\right\}}{\delta^2}=
%\frac{var (\hat{\sigma}^2_{Dicker}|{ X})+ \left\{E (\hat{\sigma}^2_{Dicker}|{ X}) - \sigma^2\right\}^2}{\delta^2}.
%\]

Now, define $\ve^{(1)}_n=\sqrt{n} {\rm E}\{{\rm var}_X(\hat{\sigma}^2_{Dicker})\}$ and  $\ve^{(2)}_n=\sqrt{n}{\rm var}\{ {\rm E}_X(\hat{\sigma}^2_{Dicker})\}$. We have that 
\[
\ve^{(1)}_n=\sqrt{n} {\rm E}\{{\rm var}_X(\hat{\sigma}^2_{Dicker})\} \le \sqrt{n} {\rm var}(\hat{\sigma}^2_{Dicker}) \le C_1 V/ \sqrt{n}, 
\]
where the first inequality follows from \eqref{eq:E_Var} and the second follows from \eqref{eq:V_bound}. Similarly, $\ve^{(2)}_n \le C_1 V/ \sqrt{n}$.  

By Markov's inequality,
\[
{\rm pr}\left\{{\rm var}_X (\hat{\sigma}^2_{Dicker}) \ge  \ve^{(1)}_n \right\} \le \frac{{\rm E}\{{\rm var}_X (\hat{\sigma}^2_{Dicker})\} }{\ve^{(1)}_n} = \frac{1}{\sqrt{n}},
\]
and also ${\rm pr}\left[ \left\{{\rm E}_X (\hat{\sigma}^2_{Dicker}) - \sigma^2\right\}^2 \ge \ve^{(2)}_n \right] \le 1/\sqrt{n}$. Therefore, for
\[
\tilde{A}_n = \left\{ { X} : {\rm var}_X (\hat{\sigma}^2_{Dicker}) \le  \ve^{(1)}_n \text{ , }
\left\{{\rm E}_X (\hat{\sigma}^2_{Dicker}) - \sigma^2\right\}^2 \le \ve^{(2)}_n \right\};
\]
we have that ${\rm pr}( { X} \in \tilde{A}_n) \to 1$. 
%We now compute the conditional probability of the event that $\delta_{Dicker}=1$ but $\sigma^2_0$ is actually true, i,e., $\sigma^2=\sigma_0^2$,
%\[
%{\rm pr}_X\{\hat{\sigma}^2_{Dicker} \ge  (\sigma^2_1+\sigma^2)/2 \} \le {\rm pr}_X\left(|\hat{\sigma}^2_{Dicker}-\sigma^2| \ge \xi \right),  
%\]
%where $\xi=(\sigma^2_1+\sigma^2)/2$. 
Now,  
for ${ X} \in \tilde{A}_n$,  the conditional Markov's inequality implies that
\begin{multline}\label{eq:tildeA}
{\rm pr}_X\left(|\hat{\sigma}^2_{Dicker}-\sigma^2| \ge \xi\right) \le \frac{{\rm E}_X\left\{\left(\hat{\sigma}^2_{Dicker}-\sigma^2\right)^2\right\}}{\xi^2}=
\frac{{\rm var}_X(\hat{\sigma}^2_{Dicker})+ \left\{{\rm E}_X(\hat{\sigma}^2_{Dicker}) - \sigma^2\right\}^2}{\xi^2}\\
\le \frac{\ve_n{(1)}+\ve_n{(2)}}{\xi^2}=\frac{2 C_1 V}{\xi^2 \sqrt{n}}.
\end{multline}

The next step is to bound $V$.
Recall the notation ${ \mu}={ X}{ \beta}$. We have that $\frac{1}{n}||{ \mu}||^2 -||{ \beta}||^2\tendp 0$. Indeed, $\frac{1}{n}{\rm E}(||{ \mu}||^2)=||{ \beta}||^2$ and 
\[
{\rm var}\left(\frac{1}{n}||{ \mu}||^2\right) =\frac{{\rm E}\left\{ \left({ x}^T { \beta}\right)^4 \right\}-||{ \beta}||^4}{n}=\frac{2||{ \beta}||^4}{n},
\]
since
\[
{\rm E}\left\{ \left({ x}^T { \beta}\right)^4 \right\}= \sum_j {\rm E}(x_j^4) \beta_j ^2 + 3\sum_{j \ne j'} \beta_j^2 \beta_{j'}^2 =  \sum_j \{{\rm E}(x_j^4) -3\} \beta_j^4 + 3||{ \beta}||^4=2||{ \beta}||^4,
\]
where in the second equality the identity $\sum_{j \ne j'} \beta_j^2 \beta_{j'}^2 + \sum_j \beta_j^4 = ||{ \beta}||^4$ is used and the last equality follows from normality. Hence, $\frac{1}{n}||{ \mu}||^2 -||{ \beta}||^2\tendp 0$ as $n \to \infty$. Therefore, also 
$\left(\frac{1}{n}||{ \mu}||^2\right)^2 -||{ \beta}||^4\tendp 0$ as $n \to \infty$. 

We have that
\begin{multline*}
V=2\left\{c(\sigma^4+||{ \beta}||^4+2\sigma^2||{ \beta}||^2)+\sigma^4+||{ \beta}||^4\right\}\\=2\left[c\left\{||{ \beta}||^4-\left(\frac{1}{n}||{ \mu}||^2\right)^2+2\sigma^2\left(||{ \beta}||^2-\frac{1}{n}||{ \mu}||^2\right)\right\}+||{ \beta}||^4 -\left(\frac{1}{n}||{ \mu}||^2\right)^2\right]\\
+2(c+1)\left\{\left(\frac{1}{n}||{ \mu}||^2\right)^2+\sigma^4\right\}+\frac{4 \sigma^2}{n}||{ \mu}||^2 + 2\sigma^4.
\end{multline*}
Since $\frac{1}{n}||{ \mu}||^2 -||{ \beta}||^2 \tendp 0$ and $\left(\frac{1}{n}||{ \mu}||^2\right)^2-||{ \beta}||^4 \tendp 0$ as $n \to \infty$, then
\begin{equation}\label{eq:barA}
\bar{A}_n=\left\{ { X} | V \le 1+ 2(c+1)\left\{\left(\frac{1}{n}||{ \mu}||^2\right)^2+\sigma^4\right\}+\frac{4 \sigma^2}{n}||{ \mu}||^2 + 2\sigma^4 \right\},
\end{equation}
satisfies ${\rm pr}( { X} \in \bar{A}_n) \to 1$ as $n \to \infty$.

Let $A_n = \tilde{A}_n \cap \bar{A}_n$, then ${\rm pr}({ X}\in A_n) \to 1$ as $n \to \infty$. It follows from \eqref{eq:tildeA} and \eqref{eq:barA}, that for ${ X} \in A_n$ 
$$
{\rm pr}_X\left(|\hat{\sigma}^2_{Dicker}-\sigma^2| \ge \xi\right) \le \frac{2 C_1 V}{\xi^2 \sqrt{n}}
\le 2 C_1 \frac{1+ 2(c+1)\left\{\left(\frac{1}{n}||{ \mu}||^2\right)^2+\sigma^4\right\}+\frac{4 \sigma^2}{n}||{ \mu}||^2 + 2\sigma^4}{\xi^2 \sqrt{n}}. %\qed
$$ 

\section{Simulations}\label{sec:sim}

The discussion in the paper was based on asymptotic arguments. Here the purpose is to evaluate $R_{\delta_{Dicker}}$ for certain values of parameters. We use a similar Bayesian setting as in Section 3, where $\mu_i$ is sampled from a mixture of normals $N(0,\eta_0^2)$ and $N(0,\eta_1^2)$, but here the $\mu$'s are not independent as indicated below. 
Table \ref{tab} reports simulation estimates of $E_{\pi'}(R_{\delta_{Dicker}})$ under various values of the parameters, where $\pi'$ is the sampling distribution of $\sigma^2$ and $\mu$. For each scenario, ${ X}$ was sampled once from a standard normal distribution with independent entries. In all the scenarios considered here, $\eta_0^2+\sigma_0^2=\eta_1^2+\sigma_1^2$; also $p=n$. For each scenario, 10,000 simulated data sets were sampled where in 5,000 data sets the ${\mu}$'s are sampled from $N(0,\eta_0^2)$ and $\sigma^2=\sigma_0^2$, while in the other 5,000  data sets, the distribution is $N(0,\eta_1^2)$ and $\sigma^2=\sigma_1^2$. Each ${ x}_i$ was standardized so that $\sum_{j} x_{ij}=0$ and $\sum_j x_{ij}^2=p$. Hence, if $\beta_j$ are i.i.d $N(0,\eta^2/p)$ then each $\mu_i$ is distributed $N(0,\eta^2)$. Note, however, that the $\mu$'s under this prior are not independent. In each data set, it was recorded whether $\delta_{Dicker}$ correctly identifies $\sigma^2$. The mean of those indicators over the simulated data sets is a simulation estimate of $E_{\pi'}(R_{\delta_{Dicker}})$.

\begin{table} [ht!]
	\caption{ Simulation estimates of $E_{\pi'}(R_{\delta_{Dicker}})$  under various settings for $\eta_1^2,\sigma_1^2$. The values of $\eta_0^2,\sigma_0^2$ are fixed to be $(1,1)$. Simulations' standard errors are given in parentheses.}
	\begin{center}
		{
			\begin{tabular}{c||c|c|c|c|c}
				%\hline
				$n=p$ & $(\eta_1^2,\sigma_1^2)=(\frac{5}{6},\frac{7}{6})$ & $(\eta_1^2,\sigma_1^2)=(\frac{4}{6},\frac{8}{6})$ & $(\eta_1^2,\sigma_1^2)=(\frac{3}{6},\frac{9}{6})$ & $(\eta_1^2,\sigma_1^2)=(\frac{1}{6},\frac{10}{6})$ &
				$(\eta_1^2,\sigma_1^2)=(\frac{1}{6},\frac{11}{6})$ \\
				\hline 
				\hline
				100&0$\cdot$396(0$\cdot$005)&0$\cdot$318(0$\cdot$005)&0$\cdot$236(0$\cdot$004)&0$\cdot$178(0$\cdot$004)&0$\cdot$123(0$\cdot$003)\\
				200&0$\cdot$363(0$\cdot$005)&0$\cdot$251(0$\cdot$004)&0$\cdot$154(0$\cdot$004)&0$\cdot$093(0$\cdot$003)&0$\cdot$046(0$\cdot$002)\\
				300&0$\cdot$338(0$\cdot$005)&0$\cdot$198(0$\cdot$004)&0$\cdot$109(0$\cdot$003)&0$\cdot$05(0$\cdot$002)&0$\cdot$02(0$\cdot$001)\\
				400&0$\cdot$31(0$\cdot$005)&0$\cdot$177(0$\cdot$004)&0$\cdot$074(0$\cdot$003)&0$\cdot$029(0$\cdot$002)&0$\cdot$008(0$\cdot$001)\\
				500&0$\cdot$293(0$\cdot$005)&0$\cdot$136(0$\cdot$003)&0$\cdot$054(0$\cdot$002)&0$\cdot$017(0$\cdot$001)&0$\cdot$004(0$\cdot$001)\\
				600&0$\cdot$279(0$\cdot$004)&0$\cdot$123(0$\cdot$003)&0$\cdot$041(0$\cdot$002)&0$\cdot$01(0$\cdot$001)&0$\cdot$002(0$\cdot$001)\\
				700&0$\cdot$267(0$\cdot$004)&0$\cdot$1(0$\cdot$003)&0$\cdot$028(0$\cdot$002)&0$\cdot$006(0$\cdot$001)&0(0)\\
				800&0$\cdot$248(0$\cdot$004)&0$\cdot$086(0$\cdot$003)&0$\cdot$021(0$\cdot$001)&0$\cdot$004(0$\cdot$001)&0(0)\\
				900&0$\cdot$238(0$\cdot$004)&0$\cdot$07(0$\cdot$003)&0$\cdot$017(0$\cdot$001)&0$\cdot$002(0$\cdot$001)&0(0)\\
				1000&0$\cdot$219(0$\cdot$004)&0$\cdot$064(0$\cdot$002)&0$\cdot$011(0$\cdot$001)&0$\cdot$002(0$\cdot$001)&0(0)
		\end{tabular}}\label{tab}
	\end{center}
\end{table}

In all scenarios $\sigma_0^2=1$ but $\sigma_1^2$ varies from $7/6$ to $11/6$. Similarly, $\eta_0^2$ is fixed to be 1 and $\eta_1^2$ goes from $5/6$ to $1/6$. As $\sigma_1^2$ becomes farther from $\sigma_0^2$ and as $n$ grows, $E_{\pi'}(R_{\delta_{Dicker}})$ gets smaller.
When $\sigma_1^2=7/6$ it is the closest to $\sigma_0^2=1$, but even then and even when $n=p=100$, $R_{\delta_{Dicker}}$ is smaller than 1/2. 
This demonstrates that even for relatively small values of $n$ the conditional error is small in this setting. However, when the $\mu$'s are sampled from $\pi$ (see Section 3), the Bayes risk of $\delta_{Dicker}$ is not smaller than the Bayes risk, which is 0$\cdot$5.
%the asymptotic findings hold in the sense that $R_{\delta_{Dicker}}$ is smaller than the Bayes risk with high probability. 

\begin{figure}[ht!]
	\begin{center}
		\includegraphics[width=.6\textwidth]{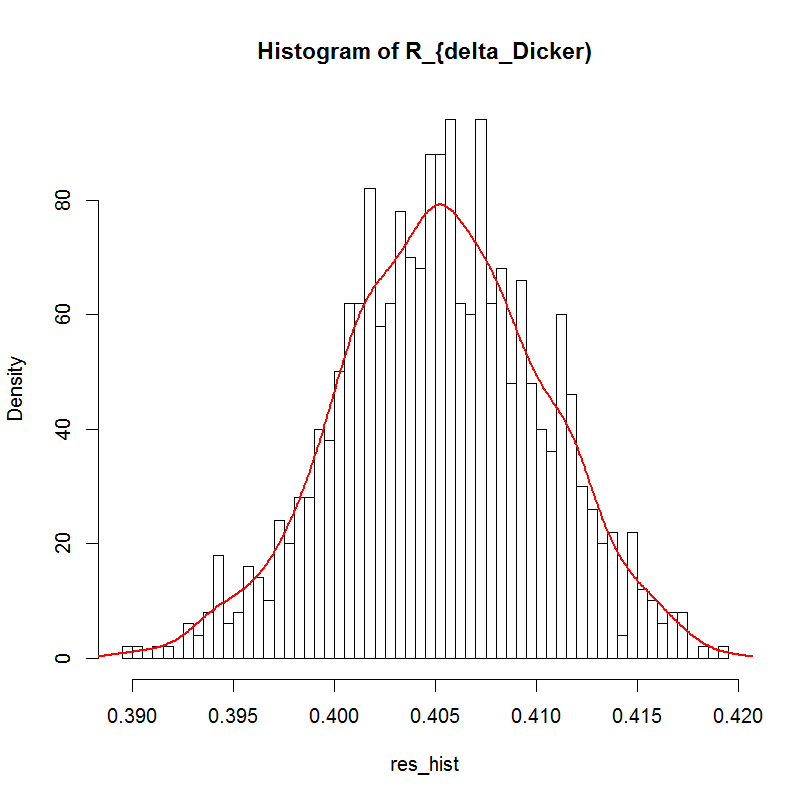}
		\caption{ Histogram of the distribution of $E_{\pi'}(R_{\delta_{Dicker}})$; the red line is a density estimate of the distribution. } \label{fig}
	\end{center}
\end{figure}

The results of Table \ref{tab} are random since they depend on the matrix ${ X}$; different ${ X}$'s yield different results for $R_{\delta_{Dicker}}$. In order to assess this randomness, the scenario $n=p=100$ with $(\eta_1^2,\sigma_1^2)=(5/6,7/6)$ was repeated 1000 times, where in each repetition 10,000 data sets were simulated as above. The resulting histogram is given in Figure \ref{fig}. It is demonstrated that in this setting, where $n$ and $p$ are not large and also  $\sigma_0^2$ and $\sigma_1^2$ are relatively close, the probability that $E_{\pi'}(R_{\delta_{Dicker}})$ is smaller than 0$\cdot$5, is practically 1.  

%Since
%\[
%{\rm E}\left(||{\bf y}^2|| \Big| {\bf X} \right)=\sum_{i=1}^n {\rm E}(y_i^2|{\bf x}_i)=\sum_{i=1}^n \left( {\bf x}_i^T { \beta}\right)^2 + n \sigma^2 = || {\bf X} { \beta}||^2 + n \sigma^2,
%\]
%and 
%\[
%{\rm E}\left(||{\bf X}^T{\bf y}||^2 \Big| {\bf X} \right)=||{ \beta}^T {\bf X}^T {\bf X}||^2+\sigma^2||{\bf X}||^2_F,
%\]
%where $||{\bf A}||_F$ is the Frobenius norm of matrix ${\bf A}$,
%then,
%\begin{multline*}
%{\rm E}(\hat{\sigma}^2_{Dicker}|{\bf X}) = \frac{p+n+1}{n(n+1)}\left(|| {\bf X} { \beta}||^2 + n \sigma^2\right)-\frac{||{ \beta}^T {\bf X}^T {\bf X}||^2+||{\bf X}||^2_F\sigma^2}{n(n+1)}\\
%=\sigma^2\frac{n(p+n+1)-||{\bf X}||^2_F}{n(n+1)}+\frac{(p+n+1)|| {\bf X} { \beta}||^2-||{ \beta}^T {\bf X}^T {\bf X}||^2 }{n(n+1)}
%\end{multline*}
%We can write 
%\[
%||{\bf X} { \beta}||^2 = { \beta}^T {\bf X}^T {\bf X} { \beta} \text{ and } ||{ \beta}^T {\bf X}^T {\bf X}||^2= { \beta}^T \left({\bf X}^T {\bf X}\right)^2 { \beta}
%\]

%\begin{thebibliography}{9}

\section*{Acknowledgment}
I am grateful to Pavel Chigansky,  Yair Goldberg and Micha Mandel who made comments on an earlier version of the manuscript. I wish to thank also the editor of Biometrika, the associate editor and two referees for their close reading and useful comments.

\bibliographystyle{apa}

\bibliography{refs}

\end{document}